\documentclass[11pt]{amsart}

\usepackage{amssymb}
\usepackage{mathrsfs}
\usepackage{amsfonts}
\usepackage{latexsym,amsmath,amsthm,amssymb,amsfonts}
\usepackage[usenames]{color}
\usepackage{xspace,colortbl}
\usepackage{graphicx}
\usepackage{tipa}
\usepackage{amsthm}

\newcommand{\be}{\begin{equation}}
\newcommand{\ee}{\end{equation}}
\newcommand{\beq}{\begin{eqnarray}}
\newcommand{\eeq}{\end{eqnarray}}

\newtheorem{prop}{Proposition}[section]

\newtheorem{remark}[prop]{Remark}

\newtheorem*{condI}{\bf Condition I}
\def\begeq{\begin{equation}}
\def\endeq{\end{equation}}

\def\odot{\setbox0=\hbox{$\bigcirc$}\relax \mathbin {\hbox
to0pt{\raise.5pt\hbox to\wd0{\hfil $\wedge$\hfil}\hss}\box0 }}

\numberwithin{equation} {section}

\numberwithin{equation}{section}
\newtheorem{theorem}{\bf Theorem}[section]
\newtheorem{proposition}[theorem]{\bf Proposition}
\newtheorem{definition}[theorem]{\bf Definition}
\newtheorem{lemma}[theorem]{\bf Lemma}

\allowdisplaybreaks

\begin{document}
\title[The Dirichlet problem for sum Hessian quotient equations on Riemannian manifolds]
 {The Dirichlet problem for sum Hessian quotient equations on Riemannian manifolds}\thanks{\it {The research is partially supported by NSFC (No. 12261105).}}


\author{
Guanghan Li,~~Chenyang Liu*}
\address{
School of Mathematics and Statistics, Wuhan University, Wuhan 430072, China
}

\email{ghli@whu.edu.cn,liuchy@whu.edu.cn}

\thanks{$\ast$ Corresponding author}

\date{}
\maketitle
\begin{abstract}
In this paper, we consider the Dirichlet problem for a class of sum Hessian quotient equations involving a gradient term on the right-hand side on Riemannian manifolds. Under the assumption of the existence of a subsolution, based on the a priori estimates for the solutions, we prove the existence and the uniqueness of $(k-1)$-admissible solutions for the Dirichlet problem on compact Riemannian manifolds with smooth boundary. Our
existence results cover the case $0\le l<k<n$ with $k\ge 2$, where the subsolution is assumed to be $(k-1)$-admissible, and the case $k=n$ with $0\le l\le n-2$, where the subsolution is assumed to be $n$-admissible.
\end{abstract}

\section{Introduction} \label{S1}

Let $(\overline{M},g)$ be a compact Riemannian manifold of dimension $n\ge 3$ with smooth boundary $\partial M$ and $\overline{M}=M\cup\partial M$. In this paper, we consider the Dirichlet problem for a class of sum Hessian quotient equations
\begin{equation}\label{Hessian quo-1}
\left\{
\begin{aligned}
&\left(\frac{\sigma_{k}+\alpha\sigma_{k-1}}{\sigma_{l}+\alpha\sigma_{l-1}}\right)^{\frac{1}{k-l}}(\eta[U])=f(x,u,\nabla u),\qquad&&~\mathrm{in}~M,\\
&u=\varphi, \qquad&&~\mathrm{on}~ \partial M,\\
\end{aligned}
\right.
\end{equation}
where $u(x)$ is a function defined over $\overline{M}$, $f$ is a positive $C^{\infty}$ function with respect to $(x, z, p)\in \overline{M}\times\mathbb{R}\times T_{x}M$, $\alpha\geq 0$, $U=(\tau\Delta u)I - \nabla^{2}u$ with $\tau\geq1$, $\eta[U]= (\eta_1,\dots,\eta_n)$ are the eigenvalues of $U$, and $\nabla$, $\Delta$, $\nabla^2$ denote the gradient, Laplace, Hessian operators on $M$ induced by $g$, respectively. 

Recently, there have been many works on equations of sum Hessian type. Ren-Wang \cite{RW25-1} considered the following sum Hessian equation:
\begin{equation}\label{Hessian-1}
\sigma_{k}(\tilde{\eta})+\alpha\sigma_{k-1}(\tilde{\eta})=f(x,u,\nabla u),
\end{equation}
where $\tilde{\eta} = (\tilde{\eta}_1,\dots,\tilde{\eta}_n)$ are the eigenvalues of $(\Delta u)I - \nabla^{2}u$, $f$ is a positive function, $\alpha\geq0$ and $0<k\leq n$. For $0< k < n$, they proved interior estimates for \eqref{Hessian-1} and Pogorelov type $C^{2}$ estimates for its Dirichlet problem on bounded domains in $\mathbb{R}^{n}$. In the case $k=n$ they obtained weak Pogorelov type $C^{2}$ estimates. In \cite{RW25-2}, they also studied the sum Hessian quotient equation and obtained similar results. When the eigenvalues are replaced by $\lambda[\nabla^{2}u]$ (i.e., the eigenvalues of $\nabla^{2}u$), Liu–Ren \cite{LR23} established Pogorelov type $C^{2}$ estimates for \eqref{Hessian-1} on bounded domains in $\mathbb{R}^{n}$ with Dirichlet boundary condition, and obtained the corresponding rigidity theorem. Dong, Xu and Zhang \cite{DXZ26} studied the general sum-type Hessian equation
\begin{equation*}
\sigma_{k}(\lambda[\nabla^{2}u])+\sum_{r=1}^{m}a_{r}\sigma_{k-r}(\lambda[\nabla^{2}u])=f(x,u,\nabla u).
\end{equation*}
They established Pogorelov type $C^{2}$ estimates under the dynamic semi-convexity assumption, and proved a Liouville-type theorem. In our previous work \cite{LL26}, we established interior $C^{2}$ estimates and Pogorelov type $C^{2}$ estimates for sum Hessian quotient equations on Riemannian manifolds. Since interior $C^{2}$ estimates play a key role in the existence theory,  a natural question is whether the existence of solutions to the sum Hessian type equations can be established.

The existence of solutions to Hessian equations and Hessian quotient equations, i.e., the case $\alpha=0$, has been extensively studied. In Euclidean space,  Ivochkina \cite{Iv83} first proved the classical solvability of the Dirichlet problem for the Monge–Amp\`ere equation $\det(\lambda[\nabla^{2}u])=f(x)$ on uniformly convex domains. Caffarelli–Nirenberg–Spruck \cite{CNS85} studied the Dirichlet problem for fully nonlinear elliptic equations of the form $F(\lambda[\nabla^{2}u])=f(x,u,\nabla u)$, and established the existence of smooth admissible solutions. Chou and Wang \cite{CW01} developed a variational theory for the Hessian equation
\begin{equation*}
\sigma_{k}(\lambda[\nabla^{2}u])=f(x,u),
\end{equation*}
and obtained the existence and regularity of $k$-convex solutions. Trudinger \cite{T90} solved the Dirichlet problem for the prescribed curvature equations, which include the Hessian quotient equations
\begin{equation*}
\Big(\frac{\sigma_{k}}{\sigma_{l}}\Big)^{\frac{1}{k-l}}(\lambda[\nabla^{2}u])=f(x).
\end{equation*}
Chen-Tu-Xiang~\cite{CTX20} proved the existence of closed star-shaped admissible hypersurfaces in $\mathbb{R}^{n+1}$ satisfying a Hessian quotient equation for the eigenvalues of
$\eta_{ij}=Hg_{ij}-h_{ij}$, where $g_{ij}$ and $h_{ij}$ are the components of the induced metric and the second fundamental form, respectively, and $H=g^{ij}h_{ij}$ is the mean curvature.
On Riemannian manifolds, Guan \cite{Guan99} first solved the Dirichlet problem for Hessian type equations $F(\lambda[\nabla^{2}u+\gamma ug])=\psi(x,u,\nabla u)$. Guan and Jiao \cite{GJ16} extended these results to Hessian equations 
\begin{equation*}
\sigma_{k}(\lambda[\nabla^{2}u+\chi(x,u,\nabla u)])=\psi(x,u,\nabla u), 
\end{equation*}
where $\chi$ is a general $(0,2)$-tensor. For Hessian quotient equations, Chen, Tu and Xiang \cite{CTX23} proved the existence and uniqueness of admissible solutions to the Dirichlet problem on Riemannian manifolds under a subsolution assumption. Guan \cite{Guan23} proved that on Riemannian manifolds, a subsolution ensures a unique smooth solution to the Dirichlet problem for a class of fully nonlinear elliptic equations.

In the last few years, researchers have also begun to study the existence of solutions to equations of sum Hessian type. Liang, Yan and Zhu \cite{LYZ26} studied the Neumann problem for elliptic and parabolic sum Hessian equations in Euclidean domains, and established a priori $C^{2}$ estimates and the existence of admissible solutions. Xiang and Xiong \cite{XX26} studied a class of sum Hessian quotient equations arising in conformal geometry on closed Riemannian manifolds, and proved the existence of admissible solutions. Ma, Wang and Zhu \cite{MWZ26} studied the prescribed sum Hessian curvature equation on closed star-shaped hypersurfaces in $\mathbb{R}^{n+1}$, and proved the existence and uniqueness of solutions. Motivated by the above developments, in this paper we focus on the Dirichlet boundary value problems for sum Hessian quotient equations on Riemannian manifolds.

We first give the definition of $k$-admissible.

\begin{definition}\label{def-adm}
Let $M$ be a Riemannian manifold and $u\in C^{2}(M)$. Denote by $\eta[U]=(\eta_{1},\dots,\eta_{n})$ the eigenvalues of $U=(\tau\Delta u)I-\nabla^{2}u$. A function $u$ is called $k$-admissible if $\eta[U]\in\Gamma_{k}$ for all $x\in M$, where $\Gamma_{k}$ is the Garding's cone
\begin{equation}\label{Gamma-1}
\Gamma_{k}=\{\eta\in\mathbb{R}^{n}|\sigma_{m}(\eta)>0, m=1,\dots, k\}.
\end{equation}
\end{definition}

As in \cite{LL26, RW25-2}, we define the sum Hessian function as
\begin{equation*}
S_{k}(\eta)=\sigma_{k}(\eta)+\alpha\sigma_{k-1}(\eta),
\end{equation*}
for $\eta=(\eta_{1},\dots,\eta_{n})\in\mathbb{R}^{n}$ and $1\leq k\leq n$. According to Li–Ren–Wang \cite{LRW19}, the sum Hessian operator $S_{k}(\eta)$ has the admissible solution set:
\begin{equation*}
\widetilde{\Gamma}_{k}=\Gamma_{k-1}\cap \{\eta|S_{k}(\eta)>0\}.
\end{equation*}
Moreover, both $S_{k}^{1/k}(\eta)$ and $[\frac{S_{k}(\eta)}{S_{l}(\eta)}]^{1/k-l}$ are concave in the set $\widetilde{\Gamma}_{k}$. 

The following theorems give the $C^{2}$ estimates for solutions of \eqref{Hessian quo-1}.

\begin{theorem}\label{maintheorem1}
Let $0\leq l<k<n$, $k\geq 2$, $\varphi\in C^{4}(\partial M)$, $f\in C^{2}(\overline{M}\times\mathbb{R}\times T_{x}M)$ with $f$, $f_{z}>0$, $u\in C^{4}(M)\cap C^{2}(\overline{M})$ be a $(k-1)$-admissible solution of the Dirichlet problem \eqref{Hessian quo-1}. Assume that there exists a $(k-1)$-admissible subsolution $\underline{u}\in C^{2}(\overline{M})$ satisfying
\begin{equation}\label{Hessian quo-2}
\left\{
\begin{aligned}
&\left(\frac{\sigma_{k}+\alpha\sigma_{k-1}}{\sigma_{l}+\alpha\sigma_{l-1}}\right)^{\frac{1}{k-l}}(\eta[\underline{U}])\geq f(x,\underline{u},\nabla \underline{u}),\qquad&&~\mathrm{in}~M,\\
&\underline{u}=\varphi, \qquad&&~\mathrm{on}~ \partial M,\\
\end{aligned}
\right.
\end{equation}
where $\underline{U}=(\tau\Delta \underline{u})I-\nabla^{2}\underline{u}$. Then we have
\begin{equation*}
\sup_{\overline{M}}|\nabla^{2}u|\leq C,
\end{equation*}
where $C$ is a positive constant depending on $n$, $k$, $l$, $\lVert u\rVert_{C^{1}}$, $\lVert\underline{u}\rVert_{C^{2}}$, $\inf f$, $\lVert f\rVert_{C^{2}}$, the metric $g$ and the manifold $M$. 
\end{theorem}

\begin{theorem}\label{maintheorem2}
Let $k=n$, $0\leq l\le n-2$, $\varphi\in C^{4}(\partial M)$, $f\in C^{2}(\overline{M}\times\mathbb{R}\times T_{x}M)$ with $f$, $f_{z}>0$, $u\in C^{4}(M)\cap C^{2}(\overline{M})$ be a $(n-1)$-admissible solution of the Dirichlet problem \eqref{Hessian quo-1}. Assume that there exists a $n$-admissible subsolution $\underline{u}\in C^{2}(\overline{M})$ satisfying \eqref{Hessian quo-2}. Then we have
\begin{equation*}
\sup_{\overline{M}}|\nabla^{2}u|\leq C,
\end{equation*}
where $C$ is a positive constant depending on $n$, $k$, $l$, $\lVert u\rVert_{C^{1}}$, $\lVert\underline{u}\rVert_{C^{2}}$, $\inf f$, $\lVert f\rVert_{C^{2}}$, the metric $g$ and the manifold $M$. 
\end{theorem}

Furthermore, to derive the gradient estimates following the method of Guan and Jiao \cite{GJ16} (see also Chen–Tu–Xiang \cite{CTX23}), we impose the following growth condition.

\begin{condI}\label{cond-I}
There exists a positive constant $K_{0}$ such that for any $x\in M$, $z\in\mathbb{R}$, $p\in T_{x}M$ with $|p|\ge K_{0}$,
\begin{equation*}
\frac{p\cdot\nabla_{x}f(x,z,p)}{|p|^{2}}-p\cdot\nabla_{p}f(x,z,p)\ge -\bar{f}(x,z)|p|^{\bar{\gamma}},
\end{equation*}
for a positive continuous function $\bar{f}$ and a constant $\bar{\gamma}\in(0,2)$.
\end{condI}

Next, we list the main existence theorems.

\begin{theorem}\label{maintheorem3}
Let $0\leq l<k<n$, $k\geq 2$, $\varphi\in C^{\infty}(\partial M)$, $f\in C^{\infty}(\overline{M}\times\mathbb{R}\times T_{x}M)$ with $f$,$f_{z}>0$. Suppose that Condition I holds and that there exists a $(k-1)$-admissible subsolution $\underline{u}\in C^{2}(\overline{M})$ satisfying \eqref{Hessian quo-2}. Then the Dirichlet problem \eqref{Hessian quo-1} has a unique $(k-1)$-admissible solution $u\in C^{\infty}(\overline{M})$.
\end{theorem}

\begin{theorem}\label{maintheorem4}
Let $k=n$, $0\leq l\le n-2$, $\varphi\in C^{\infty}(\partial M)$, $f\in C^{\infty}(\overline{M}\times\mathbb{R}\times T_{x}M)$ with $f$,$f_{z}>0$. Suppose that Condition I holds and that there exists a $n$-admissible subsolution $\underline{u}\in C^{2}(\overline{M})$ satisfying \eqref{Hessian quo-2}. Then the Dirichlet problem \eqref{Hessian quo-1} has a unique $(n-1)$-admissible solution $u\in C^{\infty}(\overline{M})$.
\end{theorem}

\begin{remark}
When $\alpha=0$, we have $S_{k}=\sigma_{k}$, so the operator in \eqref{Hessian quo-1} reduces to the classical Hessian quotient $(\sigma_k/\sigma_l)^{1/(k-l)}(\eta[U])=f$.
Moreover, the positivity of $f$ implies that $\sigma_k(\eta[U])>0$. Since $\eta[U]\in\Gamma_{k-1}$, it follows that $\eta[U]\in\Gamma_k$. Thus the solution is $k$-admissible; in particular, when $k=n$, the Dirichlet problem admits a unique $n$-admissible solution. For $0\leq l\leq k-2$, this recovers the existence and uniqueness result of Chen-Tu-Xiang~\cite{CTX23}.
\end{remark}

\begin{remark}
In particular, if $f=f(x,u)$ is independent of $p$, Condition I is not needed for the gradient estimates. Combining these with the second-order estimates yields the existence result for \eqref{Hessian quo-1}.
\end{remark}

The operator with the eigenvalue matrix $\eta[U]$ is connected with many problems in geometry. The matrix $(\Delta u)I-\nabla^{2}u$ first appeared in the study of the Gauduchon conjecture \cite{Gau84} in complex geometry, which was solved by Sz\'ekelyhidi–Tosatti–Weinkove \cite{STW17}. Harvey and Lawson \cite{HL11, HL12} introduced the notion of $(n-1)$-plurisubharmonic functions and developed the Dirichlet duality theory, where the matrix $(\Delta u)I-\nabla^{2}u$ is nonnegative definite. When $k=n$, $l=0$, $\tau=1$ and $\alpha=0$, equation \eqref{Hessian quo-1} reduces to the $(n-1)$-Monge-Amp\`ere equation
\begin{equation*}
\det\big((\Delta u)I-\nabla^{2}u\big)=f^{n}.
\end{equation*}
For strictly pseudoconvex domains in $\mathbb{C}^{n}$, Li \cite{L04} solved the Dirichlet problem, while Tosatti and Weinkove \cite{TW17} studied the equation on compact K\"ahler manifolds. For further results on related equations in the complex and real settings, one can refer to \cite{Dong23, FWW10, S18}. If the eigenvalues $\eta[U]$ are replaced by the principal curvatures $\kappa$ of a hypersurface, then \eqref{Hessian quo-1} becomes a prescribed curvature equation of quotient type, which has been widely studied; see, e.g., Caffarelli–Nirenberg–Spruck \cite{CNS85} and Guan–Ren–Wang \cite{GRW15}. For the geometric notion of $(n-1)$-convexity related to the case $k=n$, $l=0$, we refer to \cite{Sha86, Wu87}.

The paper is organized as follows. In Section \ref{S2}, some useful formulas and basic properties will be listed. In Section \ref{S3}, we establish the $C^{0}$ and $C^{1}$ estimates. Section \ref{S4} is devoted to the $C^{2}$ estimates and the proofs of Theorems \ref{maintheorem1}--\ref{maintheorem4}.

\section{Preliminaries} \label{S2}
In this section, we present some preliminaries on Riemannian manifolds and the properties of elementary symmetric polynomials.

Let $g$ be the Riemannian metric on $M$, and $\nabla$, $\Delta$, $\nabla^{2}$ denote the gradient, Laplace, and Hessian operators on $M$ induced by $g$, respectively. For a $(s, r)$-tensor field $\mathscr{T}$ on $M$, its covariant derivative $\nabla\mathscr{T}$ is a $(s, r+1)$-tensor field given by
\begin{equation*}
\begin{aligned}
&\nabla\mathscr{T}(Y^{1},\dots,Y^{s},X_{1},\dots,X_{r},X)\\
&=~\nabla_{X}\mathscr{T}(Y^{1},\dots,Y^{s},X_{1},\dots,X_{r})\\
&=~X(\mathscr{T}(Y^{1},\dots,Y^{s},X_{1},\dots,X_{r}))-\mathscr{T}(\nabla_{X}Y^{1},\dots,Y^{s},X_{1},\dots,X_{r})\\
&\quad\quad-\cdots-\mathscr{T}(Y^{1},\dots,Y^{s},X_{1},\dots,\nabla_{X}X_{r}),
\end{aligned}
\end{equation*}
where $Y^{1},\dots,Y^{s}$ are smooth $1$-forms, and $X_{1},\dots,X_{r},X$ are smooth tangent vector fields. Its components in local coordinates are denoted by
\begin{equation*}
\mathscr{T}_{k_{1}\cdots k_{r}, k_{r+1}}^{l_{1}\cdots l_{s}},
\end{equation*}
where $1\leq l_{i}, k_{j}\leq n$ with $i=1, 2, \cdots, s$ and $j=1, 2, \cdots, r+1$. In the following, a comma “,” in subscript of a given tensor means doing covariant derivatives. For simplicity, we will omit it when no confusion arises. The second covariant derivative of $\mathscr{T}$ is defined analogously:
\begin{equation*}
\nabla^{2}\mathscr{T}(Y^{1},\dots,Y^{s},X_{1},\dots,X_{r},X,Y)=(\nabla_{Y}(\nabla\mathscr{T}))(Y^{1},\dots,Y^{s},X_{1},\dots,X_{r},X)
\end{equation*}
and then its components in local coordinates are denoted by
\begin{equation*}
\mathscr{T}_{k_{1}\cdots k_{r}, k_{r+1}k_{r+2}}^{l_{1}\cdots l_{s}},
\end{equation*}
where $1\leq l_{i}, k_{j}\leq n$ with $i=1, 2, \cdots, s$ and $j=1, 2, \cdots, r+2$. Similarly, the higher order covariant derivatives of $\mathscr{T}$ are given as follows
\begin{equation*}
\nabla^{3}\mathscr{T}=\nabla(\nabla^{2}\mathscr{T}),\quad \nabla^{4}\mathscr{T}=\nabla(\nabla^{3}\mathscr{T}),\cdots,
\end{equation*}
and so on. 

For any tangent vector fields $X, Y, Z$ on $M$, the Riemannian curvature $(1, 3)$-tensor $R$ induced by the metric $g$ is defined as
\begin{equation*}
R(X, Y)Z=-\nabla_{X}\nabla_{Y}Z+\nabla_{Y}\nabla_{X}Z+\nabla_{[X, Y]}Z,
\end{equation*}
where $[\cdot,\cdot]$ denotes the Lie bracket. In a local coordinate chart $\{\xi^{i}\}_{i=1}^{n}$ of $M$, the component of the curvature tensor $R$ is defined by
\begin{equation*}
R\left(\frac{\partial}{\partial \xi^{i}},\frac{\partial}{\partial \xi^{j}}\right)\frac{\partial}{\partial \xi^{k}}=R^{h}_{kij}\frac{\partial}{\partial \xi^{h}},
\end{equation*}
and
\begin{equation*}
g\left(R\left(\frac{\partial}{\partial \xi^{k}},\frac{\partial}{\partial \xi^{l}}\right)\frac{\partial}{\partial \xi^{i}}, \frac{\partial}{\partial \xi^{j}}\right)=R_{ikl}^{h}g_{hj}=R_{ijkl},
\end{equation*}
where $g_{ij}:=g\left(\frac{\partial}{\partial \xi^{i}}, \frac{\partial}{\partial \xi^{j}}\right)$. Then, we give the well-known Ricci identity
\begin{equation}\label{ricci-ident}
\mathscr{T}_{j_{1}\cdots j_{s}, lk}^{i_{1}\cdots i_{r}}-\mathscr{T}_{j_{1}\cdots j_{s}, kl}^{i_{1}\cdots i_{r}}=-\sum\limits_{a=1}^{r}\mathscr{T}_{j_{1}\cdots j_{s}}^{i_{1}\cdots i_{a-1}ii_{a+1}\cdots i_{r}}R^{i_{a}}_{ikl}
+\sum\limits_{b=1}^{s}\mathscr{T}_{j_{1}\cdots j_{b-1}jj_{b+1}\cdots j_{s}}^{i_{1}\cdots i_{r}}R^{j}_{j_{b}kl}.
\end{equation}
If $u$ is a sufficiently smooth function on $M$, its second covariant derivative is 
\begin{equation*}
\nabla_{i}\nabla_{j}u=\frac{\partial^{2}u}{\partial\xi_{i}\partial\xi_{j}}-\Gamma_{ij}^{k}\frac{\partial u}{\partial\xi_{k}},
\end{equation*}
where $\Gamma_{ij}^{k}$ are the Christoffel symbols, $\partial/\partial \xi_{i}$ denotes the ordinary partial derivative. Throughout this paper, we write $u_{ij}$ for $\nabla_i\nabla_j u$, and similarly we write $u_{ijk}$ and $u_{ijkl}$ for higher order covariant derivatives. By the Ricci identity \eqref{ricci-ident}, one derives
\begin{equation}\label{3-eq7-1}
u_{11i}=u_{i11}+R^{m}_{1i1}u_{m}.
\end{equation}
\begin{equation}\label{3-eq7-2}
u_{11ii}=u_{ii11}+2R_{ii1}^{1}u_{11}+2R_{1i1}^{i}u_{ii}+u_{j}(R_{ii1}^{j})_{,\xi_{1}}+u_{j}(R_{1i1}^{j})_{,\xi_{i}}.
\end{equation}

Next, some algebraic identities and properties of $\sigma_k$ and $S_k$ are recalled. These can be found in \cite{CTX21, LL26, LR23, RW25-1, RW25-2}, for convenience we list them here. 

Let $\eta=(\eta_{1},\dots,\eta_{n})\in\mathbb{R}^{n}$, recall the definition of elementary symmetric function for $1\leq k\leq n$,
\begin{equation}\label{sigma-1}
\sigma_{k}(\eta)=\sum\limits_{1\leq i_{1}<\dots<i_{k}\leq n}\eta_{i_{1}}\dots\eta_{i_{k}}.
\end{equation}
We adopt the convention that $\sigma_{0}=1$ and $\sigma_{k}=0$ for $k>n$ or $k<0$. The notation $(\eta|i)=(\eta_1,\dots,\eta_{i-1},\eta_{i+1},\dots,\eta_n)$ will be used throughout. From \cite[Proposition 2.1]{CTX21}, one has the following facts:

\begin{proposition}\label{prop-1}
Let $\eta=(\eta_{1},\dots,\eta_{n})\in\mathbb{R}^{n}$ and $1\leq k\leq n$, then we have\\
(1)~$\Gamma_{1}\supset\Gamma_{2}\supset\dots\supset\Gamma_{n};$\\
(2)~$\sigma_{k-1}(\eta|i)>0$ for $\eta\in\Gamma_{k}$ and $1\leq i\leq n$;\\
(3)~$\sigma_{k}(\eta)=\sigma_{k}(\eta|i)+\eta_{i}\sigma_{k-1}(\eta|i)$ for $1\leq i\leq n$;\\
(4)~If $\eta_{1}\geq\eta_{2}\geq\dots\geq\eta_{n},$ then $\sigma_{k-1}(\eta|1)\leq\sigma_{k-1}(\eta|2)\leq\dots\leq\sigma_{k-1}(\eta|n)$ for $\eta\in\Gamma_{k}$;\\
(5)~$\sum\limits_{i=1}^{n}\sigma_{k-1}(\eta|i)=(n-k+1)\sigma_{k-1}(\eta)$.
\end{proposition}

The following Newton-MacLaurin inequality will be used as well (see, e.g., \cite{LT94, T90}).

\begin{lemma}\label{prop-2}
Let $\eta\in\mathbb{R}^{n}$. For $0\leq l< k\leq n$, $r>s\geq 0$, $k\geq r$, $l\geq s$, we have
\begin{equation}\label{N-M-ieq}
\left[\frac{\sigma_{k}(\eta)/C^{k}_{n}}{\sigma_{l}(\eta)/C^{l}_{n}}\right]^{\frac{1}{k-l}}\leq\left[\frac{\sigma_{r}(\eta)/C^{r}_{n}}{\sigma_{s}(\eta)/C^{s}_{n}}\right]^{\frac{1}{r-s}},
\quad for~\eta\in\Gamma_{k}.
\end{equation}
\end{lemma}

Some useful properties are listed below (see, e.g., \cite{RW25-1, RW25-2}).

\begin{proposition}\label{prop-3}
Let $\eta=(\eta_{1},\dots,\eta_{n})\in\mathbb{R}^{n}$ and $1\leq k\leq n$, then we have\\
(1)~$S^{ii}_{k}(\eta):=\frac{\partial S_{k}(\eta)}{\partial\eta_{i}}=\sigma_{k-1}(\eta|i)+\alpha\sigma_{k-2}(\eta|i),\quad i=1,2,\dots,n;$\\
(2)~$S^{ii,jj}_{k}(\eta):=\frac{\partial^{2} S_{k}(\eta)}{\partial\eta_{i}\partial\eta_{j}}=S_{k-2}(\eta|ij),\quad i\neq j,\quad i,j=1,2,\dots,n;$\\
(3)~$S_{k}(\eta)=\eta_{i}S_{k-1}(\eta|i)+S_{k}(\eta|i),\quad i=1,2,\dots,n;$\\
(4)~$\sum\limits_{i=1}^{n}S_{k}(\eta|i)=(n-k)S_{k}(\eta)+\alpha\sigma_{k-1}(\eta);$\\
(5)~$\sum\limits_{i=1}^{n}\eta_{i}S_{k-1}(\eta|i)=kS_{k}(\eta)-\alpha\sigma_{k-1}(\eta)$.
\end{proposition}

\begin{lemma}\label{prop-4}
We have:\\
(1)~The sets $\widetilde{\Gamma}_k$ are open and convex, and $\widetilde{\Gamma}_{1}\supset\widetilde{\Gamma}_{2}\supset\dots\supset\widetilde{\Gamma}_{n};$\\
(2)~If $\eta=(\eta_{1},\dots,\eta_{n})\in\widetilde{\Gamma}_{k}$ and $\eta_{1}\geq\eta_{2}\geq\dots\geq\eta_{n}$, then $\eta_{k-1}>0$ for $2\leq k\leq n$,
\begin{equation*}
S_{k-1}(\eta|n)\geq S_{k-1}(\eta|n-1)\geq\dots\geq S_{k-1}(\eta|1)>0,
\end{equation*}
and
\begin{equation*}
S_{k-1}(\eta|k)\geq c(n,k)S_{k-1}(\eta),
\end{equation*}
where $c(n, k)$ is a positive constant only depending on $n$ and $k$, $1\leq k\leq n$;\\
(3)~If $\eta=(\eta_{1},\dots,\eta_{n})\in\widetilde{\Gamma}_{k}$, then for any $(\xi_{1}, \dots, \xi_{n})$, we have
\begin{equation*}
\sum\limits_{i, j=1}^{n}\frac{\partial^{2}[\frac{S_{k}(\eta)}{S_{l}(\eta)}]}{\partial\eta_{i}\partial\eta_{j}}\xi_{i}\xi_{j}
\leq(1-\frac{1}{k-l})\frac{\Big[\sum_{i}\frac{\partial[\frac{S_{k}(\eta)}{S_{l}(\eta)}]}{\partial\eta_{i}}\xi_{i}\Big]^{2}}{\frac{S_{k}(\eta)}{S_{l}(\eta)}}.
\end{equation*}
\end{lemma}
\begin{proof}
Statement (1) follows directly from the definition of $S_k$. The proof of (2) can be found in \cite{R24}, and that of (3) is given in \cite{RW25-2}.
\end{proof}

\begin{lemma}\label{prop-5}
Suppose that $\eta=(\eta_{1}, \dots, \eta_{n})\in\widetilde{\Gamma}_{k}$, if $1\leq l<k\leq n$, then
\begin{equation}\label{prop-5-ieq1}
\frac{S_{k-1}(\eta)S_{l}(\eta)}{C_{n+1}^{k-1}C_{n+1}^{l}}\geq\frac{S_{k}(\eta)S_{l-1}(\eta)}{C_{n+1}^{k}C_{n+1}^{l-1}}.
\end{equation}
If $0\leq l<k$ and $0\leq q<p$, then
\begin{equation}\label{prop-5-ieq2}
\Big(\frac{S_{k}(\eta)}{S_{l}(\eta)}\Big)^{\frac{1}{k-l}}\leq\Big(\frac{S_{p}(\eta)}{S_{q}(\eta)}\Big)^{\frac{1}{p-q}},
\end{equation}
for $k\geq p$ and $l\geq q$.
\end{lemma}
\begin{proof}
The proof can be found in \cite{RW25-2}.
\end{proof}

For $\lambda=(\lambda_{1}, \dots, \lambda_{n})$, let $\eta=(\eta_1, \dots, \eta_n)$ be defined by
\begin{equation*}
\eta_{i}=\tau\sum_{j=1}^{n}\lambda_{j}-\lambda_{i},
\end{equation*}
where $\tau\geq1$. When $\tau=1$, this reduces to the setting used in \cite{RW25-2} (i.e.\ $\eta_i = \sum_{j\neq i}\lambda_j$).
In the following, we list two lemmas about $S_{k}(\eta)$ for $\eta \in \widetilde{\Gamma}_{k}$, whose proofs can be found in \cite{LL26, RW25-2}.

\begin{lemma}\label{prop-6}
Suppose that $\eta=(\eta_{1}, \dots, \eta_{n})\in\widetilde{\Gamma}_{k}$, then $S_{k}^{\frac{1}{k}}(\eta)$ and $[\frac{S_{k}(\eta)}{S_{l}(\eta)}]^{1/k-l}$ are concave. Thus, for any $(\xi_{1}, \dots, \xi_{n})$ we have
\begin{equation*}
\sum\limits_{i, j=1}^{n}\frac{\partial^{2}[\frac{S_{k}(\eta)}{S_{l}(\eta)}]}{\partial\lambda_{i}\partial\lambda_{j}}\xi_{i}\xi_{j}
\leq(1-\frac{1}{k-l})\frac{\Big[\sum_{i}\frac{\partial[\frac{S_{k}(\eta)}{S_{l}(\eta)}]}{\partial\lambda_{i}}\xi_{i}\Big]^{2}}{\frac{S_{k}(\eta)}{S_{l}(\eta)}}.
\end{equation*}
\end{lemma}

\begin{lemma}\label{prop-7}
Suppose that $\eta=(\eta_{1}, \dots, \eta_{n})\in\widetilde{\Gamma}_{k}$, and there is an ordering that $\lambda_{1}\geq\lambda_{2}\geq\dots\geq\lambda_{n}$, then we have\\
(1)~$\eta_{1}\leq\eta_{2}\leq\dots\leq\eta_{n}$, and $\eta_{n-k+2}>0$ for $2\leq k\leq n$;\\
(2)~$S_{k-1}(\eta|n-k+1)\geq c(n,k)S_{k-1}(\eta)$, $1\leq k<n$;\\
(3)~$\frac{\partial\big[\frac{S_{k}(\eta)}{S_{l}(\eta)}\big]}{\partial\eta_{1}}\geq\frac{\partial\big[\frac{S_{k}(\eta)}{S_{l}(\eta)}\big]}{\partial\eta_{2}}
\geq\dots\geq\frac{\partial\big[\frac{S_{k}(\eta)}{S_{l}(\eta)}\big]}{\partial\eta_{n}}$, $0\leq l<k\leq n$;\\
(4)~$\frac{\partial\big[\frac{S_{k}(\eta)}{S_{l}(\eta)}\big]}{\partial\lambda_{1}}\leq\frac{\partial\big[\frac{S_{k}(\eta)}{S_{l}(\eta)}\big]}{\partial\lambda_{2}}
\leq\dots\leq\frac{\partial\big[\frac{S_{k}(\eta)}{S_{l}(\eta)}\big]}{\partial\lambda_{n}}$, $0\leq l<k\leq n$;\\
(5)~$\forall~1\leq i\leq n$, $\frac{\partial\big[\frac{S_{k}(\eta)}{S_{l}(\eta)}\big]}{\partial\lambda_{i}}\geq 
c(n, k, l)\frac{1}{\tau}\sum\limits_{j=1}^{n}\frac{\partial\big[\frac{S_{k}(\eta)}{S_{l}(\eta)}\big]}{\partial\lambda_{j}}$, $0\leq l<k<n$;\\
(6)~$\sum\limits_{i=1}^{n}\frac{\partial\big[\frac{S_{k}(\eta)}{S_{l}(\eta)}\big]}{\partial\lambda_{i}}=(n\tau-1)\sum\limits_{i=1}^{n}\frac{\partial\big[\frac{S_{k}(\eta)}{S_{l}(\eta)}\big]}{\partial\eta_{i}}
\geq c(n, k, l)\tau f^{k-l-1}$ for equation \eqref{Hessian quo-1}, $0\leq l<k\leq n$.
\end{lemma}

For convenience, we set $U[u]=(\tau\Delta u)I - \nabla^{2}u$ with $\tau\geq 1$, and introduce the following notations
\begin{equation*}
F(U)=\left[\frac{S_{k}(\eta[U])}{S_{l}(\eta[U])}\right]^{\frac{1}{k-l}},\quad F^{ij}=\frac{\partial F}{\partial U_{ij}},\quad F^{ij,rs}=\frac{\partial^{2} F}{\partial U_{ij}\partial U_{rs}}.
\end{equation*}
Then equation \eqref{Hessian quo-1} is expressed in the form
\begin{equation*}
F(U)=\hat{F}(\eta[U])=f(x, u, \nabla u).
\end{equation*}
Denote
\begin{equation*}
T(\nabla^{2}u)=F((\tau\Delta u)I - \nabla^{2}u),\quad T^{ij}=\frac{\partial T}{\partial u_{ij}},\quad T^{ij,rs}=\frac{\partial^{2} T}{\partial u_{ij}\partial u_{rs}}.
\end{equation*}
Then equation \eqref{Hessian quo-1} can also be written as
\begin{equation*}
T(\nabla^{2}u)=\hat{T}(\lambda[\nabla^{2}u])=f(x, u, \nabla u),
\end{equation*}
where $\lambda[\nabla^{2}u]=(\lambda_{1},\cdots,\lambda_{n})$ are the eigenvalues of $\nabla^{2}u$. And
\begin{equation*}
T^{ij}=\frac{\partial F}{\partial U_{kl}}\cdot\frac{\partial U_{kl}}{\partial u_{ij}}=\tau\sum_{k=1}^{n}F^{kk}\delta_{ij}-F^{ij}.
\end{equation*}

By the definitions of $S_{k}$ and $f$, we have
\begin{equation}\label{Eulereq-1}
T^{ij}u_{ij}=F^{ij}U_{ij}= f-\frac{\alpha f}{k-l}\left(\frac{\sigma_{k-1}}{S_k}-\frac{\sigma_{l-1}}{S_l}\right).
\end{equation}
Differentiating $F(U)=f(x,u,\nabla u)$ with respect to $\xi_{k}$, we get
\begin{equation}\label{Eulereq-2}
T^{ij}u_{ijk}=F^{ij}U_{ijk}= f_{\xi_{k}}+f_z u_k+\sum_{i=1}^{n}f_{p_i}u_{ik},
\end{equation}
where $f_{\xi_k}=\frac{\partial f}{\partial \xi_k}$ denotes the partial derivative of $f$ with respect to its first variable.

We also need the following well-known result (see, e.g., \cite{Ball84}).

\begin{lemma}\label{prop-8}
If $W=(w_{ij})$ is a real symmetric matrix, $\eta_{i}=\eta_{i}(W)$ is one of the eigenvalues $(i=1, \dots, n)$ and $F=F(W)=\hat{F}(\eta(W))$ is a symmetric function of $\eta_{1}, \dots, \eta_{n}$, then for any real symmetric matrix $A=(a_{ij})$, we have
\begin{equation}\label{prop-7-ieq1}
\frac{\partial^{2}F}{\partial w_{ij}\partial w_{st}}a_{ij}a_{st}=\frac{\partial^{2}\hat{F}}{\partial \eta_{p}\partial \eta_{q}}a_{pp}a_{qq}
+2\sum\limits_{p<q}\frac{\frac{\partial \hat{F}}{\partial \eta_{p}}-\frac{\partial \hat{F}}{\partial \eta_{q}}}{\eta_{p}-\eta_{q}}a_{pq}^{2}.
\end{equation}
Moreover, if $\hat{F}$ is concave and $\eta_{1}\geq\eta_{2}\geq\dots\geq\eta_{n}$, we have
\begin{equation}\label{prop-7-ieq2}
\frac{\partial \hat{F}}{\partial \eta_{1}}(\eta)\leq\frac{\partial \hat{F}}{\partial \eta_{2}}(\eta)\leq\dots\leq\frac{\partial \hat{F}}{\partial \eta_{n}}(\eta).
\end{equation}
\end{lemma}
\begin{proof}
Lemma 3.2 of \cite{GLM18} (see also \cite{BAn94}) gives the proof of \eqref{prop-7-ieq1}, while Lemma 2.2 of \cite{BAn94} proves \eqref{prop-7-ieq2}.
\end{proof}

\section{$C^0$ and $C^{1}$ Estimates}  \label{S3}

In this section, following the method in \cite{CTX23, GJ16}, we establish $C^{0}$ and $C^{1}$ estimates for $(k-1)$-admissible solutions of \eqref{Hessian quo-1} for a wide range of the indices $k$ and $l$.

\begin{lemma}\label{lemma1}
Let $0\le l<k\le n$, $k\ge 2$, $\varphi\in C^{4}(\partial M)$, $f\in C^{2}(\overline{M}\times\mathbb{R}\times T_{x}M)$ with $f$, $f_{z}>0$. If $u\in C^{4}(M)\cap C^{2}(\overline{M})$ is a $(k-1)$-admissible solution of \eqref{Hessian quo-1}, and there exists a $(k-1)$-admissible subsolution $\underline{u}\in C^{2}(\overline{M})$ satisfying \eqref{Hessian quo-2}, then
\begin{equation}\label{lemma1-eq1}
\sup_{x\in\overline{M}}|u(x)|\leq C,
\end{equation}
where $C$ is a positive constant depending only on $\sup_{\partial M}\varphi$ and the subsolution $\underline{u}$.
\end{lemma}

\begin{proof}
By Definition \ref{def-adm} and $k\ge 2$, any $(k-1)$-admissible solution $u$ of \eqref{Hessian quo-1} satisfies $\eta[U]\in\Gamma_{k-1}\subset\Gamma_{1}$, which gives $\operatorname{tr}(\eta[U])=(n\tau-1)\Delta u>0$. Since $\Delta u>0$, the maximum principle gives
\begin{equation*}
\sup_{\overline{M}}u\leq\sup_{\partial M}\varphi.
\end{equation*}
On the other hand, by assumption, there exists a $(k-1)$-admissible subsolution $\underline{u}\in C^{2}(\overline{M})$ satisfying \eqref{Hessian quo-2}. Since $f$ is strictly increasing in the $u$-variable ($f_{z}>0$) and $(S_{k}/S_{l})^{1/(k-l)}$ is monotone increasing in $\eta$, the comparison principle gives
\begin{equation*}
u\ge\underline{u},\qquad x\in \overline{M}.
\end{equation*}
\end{proof}

\begin{lemma}\label{lemma2}
Let $0\le l<k\le n$, $k\ge 2$, $\varphi\in C^{4}(\partial M)$, $f\in C^{2}(\overline{M}\times\mathbb{R}\times T_{x}M)$ with $f$, $f_{z}>0$. If $u\in C^{4}(M)\cap C^{2}(\overline{M})$ is a $(k-1)$-admissible solution of \eqref{Hessian quo-1} and Condition I holds, then
\begin{equation}\label{lemma2-eq1}
\sup_{M}|\nabla u|\leq C(1+\sup_{\partial M}|\nabla u|),
\end{equation}
where $C$ is a positive constant depending only on $n$, $k$, $l$, $\lVert u\rVert_{C^{0}}$, $\lVert f\rVert_{C^{1}}$, the metric $g$ and the manifold $M$.
\end{lemma}

\begin{proof}
Consider the following test function
\begin{equation*}
P(x)=\log |\nabla u|-a\log \omega,
\end{equation*}
where $\omega(u)=-u+\sup_{M}u+1$, $a\in(0, 1)$ is a constant. Suppose that $P(x)$ attains its maximum value at an interior point $x_{0}\in M$, otherwise, the desired estimate follows. Choosing a local orthonormal frame field $\{\xi_{i}\}_{i=1}^{n}$ at $x_{0}$ such that $g_{ij}(x_{0})=\delta_{ij}$, and then rotating the coordinate axes, we can diagonalize the Hessian matrix $\nabla^{2}u=(u_{ij})$.

Differentiating $P$ at $x_{0}$ once yields
\begin{equation}\label{lemma2-eq2}
0=\frac{u_{i}u_{ii}}{|\nabla u|^{2}}+\frac{a u_{i}}{\omega},
\end{equation}
and differentiating $P$ at $x_{0}$ twice,
\begin{equation}\label{lemma2-eq3}
0\geq\frac{\sum_{k}u_{k}u_{kii}+u_{ii}^{2}}{|\nabla u|^{2}}-\frac{u_{i}^{2}u_{ii}^{2}}{|\nabla u|^{4}}+\frac{(a-a^{2}) u_{i}^{2}}{\omega^{2}}+\frac{a u_{ii}}{\omega}.
\end{equation}
Assume that $|\nabla u(x_{0})|\ge c_{0}$ for some positive constant $c_{0}$, since otherwise we are done. Set
\begin{equation*}
S=\{i\in\{1,\ldots,n\}~|~u_{i}\neq0\}.
\end{equation*}
Clearly, $S\ne\emptyset$, and for each $i\in S$, equation \eqref{lemma2-eq2} yields
\begin{equation*}
u_{ii}=-\frac{a|\nabla u|^{2}}{\omega}<0.
\end{equation*}
Since $\Delta u>0$, we have $T^{ii}\geq \sum_{l\neq i}F^{ll}\geq \frac{1}{2}\sum_{l}F^{ll}$ for each $i\in S$. Thus
\begin{equation}\label{lemma2-ieq4}
T^{ii}u_{i}^{2}=\sum_{i\in S}T^{ii}u_{i}^{2}\geq\sum_{i\in S}\left(\frac{1}{2}\sum_{l}F^{ll}\right)u_{i}^{2}=\frac{1}{2(n\tau-1)}\left(\sum_{l}T^{ll}\right)|\nabla u|^{2}.
\end{equation}
Since $f_{z}>0$ and $T^{ii}>0$, multiplying both sides of \eqref{lemma2-eq3} by $T^{ii}$ and substituting \eqref{3-eq7-1}, \eqref{Eulereq-2}, \eqref{lemma2-ieq4}, and Condition I into it, we obtain
\begin{equation}\label{lemma2-ieq5}
\begin{aligned}
0&\geq \frac{(a-a^{2})T^{ii}u_{i}^{2}}{\omega^{2}}+\frac{\sum_{k}u_{k}f_{\xi_{k}}}{|\nabla u|^{2}}+f_{z}-\frac{a}{\omega}\sum_{k}u_{k}f_{p_{k}}-Ca^{2}\sum_{i}T^{ii}+\frac{a T^{ii}u_{ii}}{\omega}\\
&\geq \frac{c_{1}|\nabla u|^{2}}{2(n\tau-1)}\sum_{l}T^{ll}-C_{1}|\nabla u|^{\bar{\gamma}}-C_{1}|\nabla u|-C_{1}\sum_{i}T^{ii}+\frac{a T^{ii}u_{ii}}{\omega},
\end{aligned}
\end{equation} 
where $c_1$ and $C_1$ are positive constants depending on $\lVert u\rVert_{C^{0}}$, $\lVert f\rVert_{C^{1}}$, the metric $g$ and the manifold $M$.

To estimate the term involving $T^{ii}u_{ii}$, we define, for $(\eta,\alpha)\in\Gamma_k^{n+1}$,
\begin{equation*}
\tilde{F}(\eta,\alpha)=\left(\frac{\sigma_k(\eta,\alpha)}{\sigma_l(\eta,\alpha)}\right)^{\frac{1}{k-l}}
=\left(\frac{\sigma_{k}(\eta)+\alpha\sigma_{k-1}(\eta)}{\sigma_{l}(\eta)+\alpha\sigma_{l-1}(\eta)}\right)^{\frac{1}{k-l}},
\end{equation*}
where $\Gamma_k^{n+1}$ is the Gårding cone of $\sigma_k$ in $\mathbb{R}^{n+1}$. Note that $(\eta[U],\alpha)\in\Gamma_k^{n+1}$ and $\tilde{F}(\eta[U],\alpha)=F(U)$. Since $k\le n$, we have $(1,\ldots,1,0)\in\Gamma_k^{n+1}$, and hence we can choose a sufficiently small constant $\varepsilon=\varepsilon(n,k)>0$ such that $\mathbf{e}_\varepsilon=(1,\ldots,1,-\varepsilon)\in\Gamma_k^{n+1}$.
By the theory of Garding cones, $\tilde{F}$ is concave and homogeneous of degree one. Therefore, in an orthonormal frame diagonalizing $U$ we have
\begin{equation*}
0<\tilde F(\mathbf{e}_\varepsilon)\le \tilde F(\zeta)+D\tilde F(\zeta)\cdot(\mathbf{e}_\varepsilon-\zeta)
=D\tilde F(\zeta)\cdot\mathbf{e}_\varepsilon=\sum_iF^{ii}-\varepsilon\tilde F_\alpha(\zeta),
\end{equation*}
where $\zeta=(\eta[U],\alpha)$ and $\tilde{F}_{\alpha}=\frac{\partial\tilde{F}}{\partial\alpha}$. Hence $\tilde{F}_{\alpha}(\zeta)\le\frac{1}{\varepsilon}\sum_{i}F^{ii}$. Since $0<a<1$ and $\omega\ge1$, using \eqref{Eulereq-1}, $f>0$, and $\sum_{i}T^{ii}=(n\tau-1)\sum_{i}F^{ii}$, it follows that
\begin{equation}\label{lemma2-ieq6}
\frac{aT^{ii}u_{ii}}{\omega}=\frac{a}{\omega}(f-\alpha\tilde{F}_{\alpha})\ge -\frac{\alpha}{\varepsilon(n\tau-1)}\sum_{i}T^{ii}\ge -C\sum_{i}T^{ii}.
\end{equation}
By \eqref{lemma2-ieq5} and \eqref{lemma2-ieq6}, we have
\begin{equation}\label{lemma2-ieq7}
0\geq \left(\frac{c_{1}}{2(n\tau-1)}|\nabla u|^{2}-C_{2}\right)\sum_{l}T^{ll}-C_{1}|\nabla u|^{\bar{\gamma}}-C_{1}|\nabla u|.
\end{equation} 
If $\frac{c_{1}}{2(n\tau-1)}|\nabla u|^{2}-C_{2}\leq 0$, then $|\nabla u|\leq C$. Otherwise, by Lemma \ref{prop-7} (6) and \eqref{lemma2-ieq7}, we obtain
\begin{equation*}
0\geq c(n,k,l)\tau \left(\frac{c_{1}}{2(n\tau-1)}|\nabla u|^{2}-C_{2}\right)-C_{1}|\nabla u|^{\bar{\gamma}}-C_{1}|\nabla u|.
\end{equation*}
Since $\bar{\gamma}\in(0,2)$, it follows that $|\nabla u|\le C$, completing the proof.
\end{proof}

Next, we give the global $C^{1}$ estimates for the solution of equation \eqref{Hessian quo-1}.

\begin{lemma}\label{lemma3}
Let $0\le l<k\le n$, $k\ge 2$, $\varphi\in C^{4}(\partial M)$, $f\in C^{2}(\overline{M}\times\mathbb{R}\times T_{x}M)$ with $f$, $f_{z}>0$. If $u\in C^{4}(M)\cap C^{2}(\overline{M})$ is a $(k-1)$-admissible solution of \eqref{Hessian quo-1} and Condition I holds, then
\begin{equation}\label{lemma3-eq1}
\sup_{\overline{M}}|\nabla u|\leq C,
\end{equation}
where $C$ is a positive constant depending only on $n$, $k$, $l$, $\lVert u\rVert_{C^{0}}$, $\lVert\underline{u}\rVert_{C^{1}}$, $\lVert f\rVert_{C^{1}}$, $\lVert\varphi\rVert_{C^{1}}$, the metric $g$ and the manifold $M$.
\end{lemma}
\begin{proof}
By Lemma \ref{lemma2}, it suffices to bound $|\nabla u|$ on $\partial M$. Let $h$ be the harmonic function in $M$ with $h=\varphi$ on $\partial M$. Since $\Delta h=0$, it follows that $\Delta(u-h)>0$ in $M$, and $u-h=0$ on $\partial M$. The maximum principle then yields $u\le h$ in $\overline{M}$, and combining this with $u\ge\underline{u}$ from Lemma \ref{lemma1}, we obtain
\begin{equation*}
\underline{u}\le u\le h,\qquad x\in\overline{M}.
\end{equation*}
Since $u$, $\underline{u}$ and $h$ are all equal to $\varphi$ on $\partial M$, we have
\begin{equation*}
\nabla_{\nu}h\le\nabla_{\nu}u\le\nabla_{\nu}\underline{u}\qquad\text{on } \partial M,
\end{equation*}
where $\nu$ is the unit exterior normal to $\partial M$. Since the tangential gradient of $u$ on $\partial M$ is determined by $\varphi$, we have
\begin{equation*}
\sup_{\partial M}|\nabla u|\le C.
\end{equation*}
Combining this with Lemma \ref{lemma2}, we obtain $\sup_{\overline{M}}|\nabla u|\le C$, completing the proof.
\end{proof}

\section{$C^{2}$ Estimates}  \label{S4}

In this section, following the method in \cite{CTX23}, we prove the global $C^{2}$ estimates and give the proofs of Theorems \ref{maintheorem1} and \ref{maintheorem2}. Since the admissible set of the solution $u$ considered here is different from that of the Hessian quotient equation, for the boundary $C^{2}$ estimates we cannot derive an important lemma on a general cone as in \cite[Lemma 4.2]{CTX23}. Instead, we obtain a similar result on the set $\widetilde{\Gamma}_{k}$. We first present the interior $C^{2}$ estimate.

\begin{lemma}\label{lemma4}
Let $0\le l<k<n$ with $k\ge 2$, or $k=n$ with $0\le l\le n-2$. Let $\varphi\in C^{4}(\partial M)$, $f\in C^{2}(\overline{M}\times\mathbb{R}\times T_{x}M)$ with $f$, $f_{z}>0$, and $u\in C^{4}(M)\cap C^{2}(\overline{M})$ be a $(k-1)$-admissible solution of \eqref{Hessian quo-1}. Assume that there exists a $(k-1)$-admissible subsolution $\underline{u}\in C^{2}(\overline{M})$ satisfying \eqref{Hessian quo-2} when $k<n$, and an $n$-admissible subsolution $\underline{u}\in C^{2}(\overline{M})$ satisfying \eqref{Hessian quo-2} when $k=n$. Then
\begin{equation*}
\sup_{M}|\nabla^{2}u|\le C\Big(1+\sup_{\partial M}|\nabla^{2}u|\Big),
\end{equation*}
where $C$ is a positive constant depending only on $n$, $k$, $l$, $\lVert u\rVert_{C^{1}}$, $\lVert\underline{u}\rVert_{C^{2}}$, $\lVert f\rVert_{C^{2}}$, the metric $g$ and the manifold $M$.
\end{lemma}

\begin{proof}
Consider the following test function
\begin{equation*}
G=\log\lambda_{\max}[\nabla^{2}u]+\frac{a}{2}|\nabla u|^{2}+A(\underline{u}-u),
\end{equation*}
for $0<a\ll 1\ll A$. Suppose that $G(x)$ attains its maximum value at an interior point $x_{0}\in M$, otherwise, the desired estimate follows. Choosing a local orthonormal frame field $\{\xi_{i}\}_{i=1}^{n}$ at $x_{0}$ such that $g_{ij}(x_{0})=\delta_{ij}$, and then rotating the coordinate axes, we can diagonalize the Hessian matrix $\nabla^{2}u=(u_{ij})$.
Differentiating $G$ at $x_{0}$ once yields
\begin{equation*}
0=G_{i}=\frac{u_{11i}}{u_{11}}+a\sum_{k}u_{k}u_{ki}+A(\underline{u}_{i}-u_{i}),
\end{equation*}
and differentiating $G$ at $x_{0}$ twice,
\begin{equation*}
0\geq G_{ii}=\frac{u_{11ii}}{u_{11}}-\frac{u_{11i}^{2}}{u_{11}^{2}}+a\sum_{k}u_{k}u_{kii}+au_{ii}^{2}+A(\underline{u}_{ii}-u_{ii}).
\end{equation*}
Let $\mathcal{L}=T^{ij}\nabla_{ij}-f_{p_{i}}\nabla_{i}$, we have
\begin{equation}\label{lemma4-eq1}
\begin{aligned}
0\geq \mathcal{L}G
=&~\frac{T^{ii}u_{11ii}}{u_{11}}-\frac{T^{ii}u_{11i}^{2}}{u_{11}^{2}}+a\sum_{k}T^{ii}u_{kii}u_{k}+aT^{ii}u_{ii}^{2}\\
&+AT^{ii}(\underline{u}_{ii}-u_{ii})-f_{p_{i}}\left(\frac{u_{11i}}{u_{11}}+a\sum_{k}u_{k}u_{ki}+A(\underline{u}_{i}-u_{i})\right)
\end{aligned}
\end{equation}
We now estimate the term $AT^{ii}(\underline{u}_{ii}-u_{ii})-A\sum_{i}f_{p_{i}}(\underline{u}_{i}-u_{i})$. It follows from \eqref{Eulereq-1} and \eqref{lemma2-ieq6} that
\begin{equation}\label{lemma4-eq2}
AT^{ii}(\underline{u}_{ii}-u_{ii})=AT^{ii}\underline{u}_{ii}-A\left(f-\alpha\tilde{F}_{\alpha}\right),
\end{equation}
where $\tilde{F}_{\alpha}$ is defined in the proof of Lemma \ref{lemma2}. By \eqref{N-M-ieq} and \eqref{Eulereq-1}, we have
\begin{equation*}
\tilde{F}_{\alpha}=\frac{f}{k-l}\left(\frac{\sigma_{k-1}}{S_k}-\frac{\sigma_{l-1}}{S_l}\right)>0.
\end{equation*}
When $k<n$, we derive
\begin{equation}\label{lemma4-eq3}
AT^{ii}(\underline{u}_{ii}-u_{ii})-A\sum_{i}f_{p_{i}}(\underline{u}_{i}-u_{i})\geq-CA\sum_{i=1}^{n}T^{ii}-CA.
\end{equation}
where $C$ is a positive constant depending on $n$, $||u||_{C^{1}}$, $||\underline{u}||_{C^{2}}$, $||f||_{C^{1}}$. When $k=n$, $\underline{u}$ is $n$-admissible, so $\eta[\underline{U}]\in\Gamma_{n}$ and $\{\eta[\underline{U}](x):x\in\overline{M}\}$ is a compact subset of $\Gamma_{n}$. Hence there exists a positive constant $c_{0}$ such that $\eta_{i}[\underline{U}]\ge c_{0}$ for $1\le i\le n$. Combining this with \eqref{lemma4-eq2}, we obtain
\begin{equation}\label{lemma4-eq4}
AT^{ii}(\underline{u}_{ii}-u_{ii})-A\sum_{i}f_{p_{i}}(\underline{u}_{i}-u_{i})\geq \frac{c_0A}{n\tau-1}\sum_iT^{ii}-CA.
\end{equation}

Combining \eqref{lemma4-eq1}, \eqref{lemma4-eq3} and \eqref{lemma4-eq4}, the subsequent computation is similar to the proofs of Theorems 1.2 and 1.3 in our previous work \cite{LL26}, so we omit the details here. Thus, we can control the largest eigenvalue of $\nabla^{2}u$, completing the proof.
\end{proof}

Next, we consider the boundary $C^{2}$ estimates. For any fixed $x_{0}\in\partial M$, we choose smooth orthonormal local frames $e_{1},\dots,e_{n}$ around $x_{0}$ such that when restricted to $\partial M$, $e_{n}$ is the unit interior normal of $\partial M$. For $x\in \overline{M}$, let $\rho(x)$ and $d(x)$ denote the distances from
$x$ to $x_{0}$ and $\partial M$, respectively,
\begin{equation*}
\rho(x)=\operatorname{dist}_{M}(x,x_{0}),\qquad d(x)=\operatorname{dist}_{M}(x,\partial M),
\end{equation*}
and set $M_{\delta}=\{x\in M:\rho(x)<\delta\}$. Since $\nabla_{ij}\rho^{2}(x_{0})=2\delta_{ij}$, we may assume that $\rho$ is smooth in $M_{\delta_{0}}$ and
\begin{equation*}
I\le\nabla_{ij}\rho^{2}\le 3I\qquad\text{in}~M_{\delta_{0}},
\end{equation*}
for $\delta_{0}>0$ sufficiently small. Then we have the following important lemma, which plays a key role in our boundary estimates.

\begin{lemma}\label{role-lemma}
Under the assumptions of Lemma \ref{lemma4}, let
\begin{equation*}
\mathcal{L}=T^{ij}\nabla_{ij}-f_{p_i}\nabla_i,\qquad v=u-\underline u+td-\frac{N}{2}d^2,
\end{equation*}
then for a positive constant $\varepsilon_{0}$, there exist uniform constants $t$, $\delta$ sufficiently small and $N$ sufficiently large such that 
\begin{equation*}
\left\{
\begin{aligned}
&\mathcal{L} v\leq-\varepsilon_{0}\left(1+\sum_{i}F^{ii}\right) && \text{in}~M_{\delta},\\
&v\geq0 && \text{on}~\partial M_{\delta}.
\end{aligned}
\right.
\end{equation*}
\end{lemma}

\begin{proof}
Since $\partial M$ is smooth and compact, we can choose a constant $\delta_{0}>0$ sufficiently small such that $d$ is smooth in
\begin{equation*}
\Omega_{\delta_{0}}=\{x\in\overline{M}:d(x)\le\delta_{0}\},
\end{equation*}
where $d(x)$ is defined as before. Then on $\Omega_{\delta_0}$, we have
\begin{equation*}
|\nabla d|=1,\qquad |\nabla^{2}d|\le C,
\end{equation*}
where $C$ is a positive constant depending only on the metric $g$ and the manifold $M$.

By the admissibility of $\underline{U}$ and the compactness of $\Omega_{\delta_{0}}$, we may choose a constant $\varepsilon>0$ sufficiently small such that $\eta[\underline{U}-\varepsilon g]\in\widetilde{\Gamma}_{k}$ in $\Omega_{\delta_{0}}$. 
When $k=n$, $\underline{u}$ is $n$-admissible, so $\eta[\underline{U}(x)]\in\Gamma_{n}$ for every $x\in\Omega_{\delta_{0}}$. Since $\underline{u}\in C^{2}(\overline{M})$, the smallest eigenvalue of $\underline{U}$ is continuous, and by the compactness of $\Omega_{\delta_{0}}$ we have
\begin{equation*}
m:=\min_{x\in\Omega_{\delta_{0}}}\lambda_{\min}\big(\underline{U}(x)\big)>0.
\end{equation*}
Thus $\underline{U}\ge mg$ in $\Omega_{\delta_{0}}$. In this case we take $\varepsilon=c_{0}=m/2$, which gives $\eta[\underline{U}-\varepsilon g]\in\Gamma_{n}$ and
\begin{equation}\label{RLemma-ieq1}
\underline{U}-\varepsilon g\ge(m-\varepsilon)g=c_{0}g \qquad\text{in } \Omega_{\delta_{0}},
\end{equation}
For $N\ge0$, set
\begin{equation*}
R_{N}=\underline{U}-\varepsilon g+N\big(\tau g-\nabla d\otimes\nabla d\big).
\end{equation*}
Since $\tau\ge1$ and $|\nabla d|=1$, we have
\begin{equation*}
R_{N}\ge\underline{U}-\varepsilon g+N\big(g-\nabla d\otimes\nabla d\big)\ge\underline{U}-\varepsilon g.
\end{equation*}
Note that $\eta\in\widetilde{\Gamma}_{k}$ if and only if $(\eta,\alpha)\in\Gamma_{k}^{n+1}$, and the Garding cone is monotone in each coordinate. Thus $\eta[R_{N}]\in\widetilde{\Gamma}_{k}$ and
$\eta[\underline{U}-\varepsilon g+N(g-\nabla d\otimes\nabla d)]\in\widetilde{\Gamma}_{k}$, and the ellipticity of $F$ gives
\begin{equation}\label{RLemma-ieq2}
F(R_{N})\ge F\big(\underline{U}-\varepsilon g+N(g-\nabla d\otimes\nabla d)\big).
\end{equation}
For fixed $x\in\Omega_{\delta_{0}}$, we take $e_{n}=\nabla d(x)$ and complete it to an orthonormal basis $\{e_{1},\dots,e_{n}\}$ of $T_{x}M$, in which $g-\nabla d\otimes\nabla d$ is diagonal with entries $(1,\dots,1,0)$.
Hence, for $1\le j\le n-1$, 
\begin{equation*}
S_{j}\big(\eta[\underline{U}-\varepsilon g+N(g-\nabla d\otimes\nabla d)]\big)=C_{n-1}^{j}N^{j}+O(N^{j-1})\qquad\text{on}~\Omega_{\delta_{0}}.
\end{equation*}
For $j=0$, we have $S_{0}=1$. 
When $k<n$, we obtain
\begin{equation*}
\frac{S_{k}}{S_{l}}\left(\eta[\underline{U}-\varepsilon g+N(g-\nabla d\otimes\nabla d)]\right)=\frac{C_{n-1}^{k}N^{k}+O(N^{k-1})}{C_{n-1}^{l}N^{l}+O(N^{l-1})},
\end{equation*}
in the case $l=0$, $S_{0}=1$. Combining the above equality with \eqref{RLemma-ieq2} and $F=\left(S_{k}/S_{l}\right)^{1/(k-l)}$, we obtain $F(R_{N})\ge cN$ for all sufficiently large $N$.
When $k=n$, we have
\begin{equation*}
S_{n}\big(\eta[\underline{U}-\varepsilon g+N(g-\nabla d\otimes\nabla d)]\big)=\big((\underline{U}_{nn}-\varepsilon g_{nn})+\alpha\big)N^{n-1}+O(N^{n-2}).
\end{equation*}
Hence, by \eqref{RLemma-ieq1},
\begin{equation*}
\frac{S_{n}}{S_{l}}\big(\eta[\underline{U}-\varepsilon g+N(g-\nabla d\otimes\nabla d)]\big)\ge\frac{(c_{0}+\alpha)N^{n-1}+O(N^{n-2})}{C_{n-1}^{l}N^{l}+O(N^{l-1})}.
\end{equation*}
Therefore, for sufficiently large $N$,
\begin{equation*}
F(R_{N})\ge cN^{\frac{n-1-l}{n-l}}.
\end{equation*}
Since $l\le n-2$, this exponent is positive. Thus in both cases $F(R_{N})$ tends uniformly to infinity on $\Omega_{\delta_{0}}$.
By the concavity of $F$ on $\widetilde{\Gamma}_{k}$, we have
\begin{equation*}
F(R_{N})-f(x,u,\nabla u)\le F^{ij}(R_{N}-U)_{ij}.
\end{equation*}
where $(R_{N}-U)_{ij}$ denotes the $(i,j)$-component of the tensor $R_{N}-U$. Using $|\nabla d|=1$ and the definition of $T^{ij}$, we obtain
\begin{equation*}
F^{ij}(\underline{U}-U)_{ij}=T^{ij}(\underline{u}_{ij}-u_{ij}),\qquad F^{ij}(\tau\delta_{ij}-d_{i}d_{j})=T^{ij}d_{i}d_{j}.
\end{equation*}
It follows that
\begin{equation}\label{RLemma-eq3}
T^{ij}(u_{ij}-\underline{u}_{ij})-NT^{ij}d_{i}d_{j} \le f(x,u,\nabla u)-F(R_{N})-\varepsilon\sum_{i}F^{ii}.
\end{equation}

We now take $\delta>0$ sufficiently small such that $M_{\delta}\subset\Omega_{\delta_{0}}$, so that at each point of $M_{\delta}$ we choose the orthonormal frame as above. A direct computation gives
\begin{equation}\label{RLemma-eq4}
\mathcal{L} v=T^{ij}(u_{ij}-\underline{u}_{ij})-NT^{ij}d_{i}d_{j}+(t-Nd)T^{ij}d_{ij}-f_{p_{i}}\left((u-\underline{u})_{i}+(t-Nd)d_{i}\right).
\end{equation}
Let $\bar{f}=\sup_{\overline{M}}f(x,u,\nabla u)$, which is bounded on $\overline{M}$. Since $d\le\delta$ in $M_{\delta}$ and $\sum_{i}T^{ii}=(n\tau-1)\sum_{i}F^{ii}$, it follows from \eqref{RLemma-eq3} and $|\nabla^{2}d|\le C$ that
\begin{equation*}
\mathcal{L} v\le-\big(\varepsilon-C(t+N\delta)\big)\sum_{i}F^{ii}+\bar{f}-F(R_{N})+C(1+t+N\delta),
\end{equation*}
where $C$ is a positive constant independent of $t$, $N$ and
$\delta$.
First choose $0<t\le1$ such that $Ct\le\varepsilon/4$. Next choose $N$ sufficiently large such that $F(R_{N})\ge \bar{f}+3C+1$ on $M_{\delta}$. Finally, shrink
$\delta$ such that $M_{\delta}\subset\Omega_{\delta_{0}}$ and $N\delta\le t$. We obtain
\begin{equation*}
\mathcal{L} v\leq-\varepsilon_{0}\left(1+\sum_{i}F^{ii}\right)\qquad\text{in } M_{\delta},
\end{equation*}
with $\varepsilon_{0}=\min\{1,\varepsilon/2\}$.
Finally, by Lemma \ref{lemma1} we have $u\ge\underline{u}$ in $\overline{M}$, and $N\delta\le t$ implies
\begin{equation*}
v=u-\underline{u}+d\left(t-\frac{N}{2}d\right)\ge u-\underline{u}+\frac{t}{2}d\ge0\qquad\text{in } \overline{M}_{\delta}.
\end{equation*}
In particular, $v=0$ on $\partial M\cap B_{\delta}(x_{0})$, completing the proof.
\end{proof}

\begin{lemma}\label{lemma5}
Let $0\le l<k<n$ with $k\ge 2$, or $k=n$ with $0\le l\le n-2$. Let $\varphi\in C^{4}(\partial M)$, $f\in C^{2}(\overline{M}\times\mathbb{R}\times T_{x}M)$ with $f$, $f_{z}>0$, and $u\in C^{4}(M)\cap C^{2}(\overline{M})$ be a $(k-1)$-admissible solution of \eqref{Hessian quo-1}. Assume that there exists a $(k-1)$-admissible subsolution $\underline{u}\in C^{2}(\overline{M})$ satisfying \eqref{Hessian quo-2} when $k<n$, and an $n$-admissible subsolution $\underline{u}\in C^{2}(\overline{M})$ satisfying \eqref{Hessian quo-2} when $k=n$. Then
\begin{equation*}
\sup_{\overline{M}}|\nabla^{2}u|\le C,
\end{equation*}
where $C$ is a positive constant depending only on $n$, $k$, $l$, $\lVert u\rVert_{C^{1}}$, $\lVert\underline{u}\rVert_{C^{2}}$, $\inf f$, $\lVert f\rVert_{C^{2}}$, the metric $g$ and the manifold $M$.
\end{lemma}

\begin{proof}
By Lemma~\ref{lemma4}, we only need to derive boundary estimates. For any $x_{0}\in\partial M$, we can choose the local orthonormal frames $e_{1},\dots,e_{n}$ around $x_{0}$ such that $e_{n}$ is interior normal to $\partial M$.

\textbf{Case 1}: Estimates of $\nabla_{\beta\gamma}u$, $\beta,\gamma=1,\dots,n-1$ on $\partial M$.

Consider a point $x_{0}\in\partial M$. Since $u-\underline{u}=0$ on $\partial M$, we have
\begin{equation*}
\nabla_{\beta\gamma}(u-\underline{u})=-\nabla_{n}(u-\underline{u})B_{\beta\gamma}\qquad\text{on } \partial M,
\end{equation*}
where
$B_{\beta\gamma}=\langle\nabla_{e_{\beta}}e_{\gamma},e_{n}\rangle$ denotes the second fundamental form of $\partial M$. Therefore,
\begin{equation*}
|\nabla_{\beta\gamma}u|\le C\qquad\text{on } \partial M,
\end{equation*}
where $C$ is a positive constant depending on $\lVert u\rVert_{C^{1}}$ and $\lVert\underline{u}\rVert_{C^{2}}$.

\textbf{Case 2}: Estimates of $\nabla_{\beta n}u$, $\beta=1,\dots,n-1$, on $\partial M$. Let
\begin{equation}\label{5-eq1}
\Phi=A_1v+A_2\rho^2-A_3\sum_{\gamma<n}\bigl|\nabla_\gamma(u-\widetilde{\varphi})\bigr|^2,
\end{equation}
where $\widetilde{\varphi}$ is a fixed $C^4$ extension of $\varphi$ to a neighborhood of $\partial M$. Using Lemma~\ref{role-lemma}, we claim that
\begin{equation*}
\left\{
\begin{aligned}
&\mathcal L\bigl(\Phi\pm\nabla_\beta(u-\widetilde{\varphi})\bigr)\leq0 &&\text{in }M_\delta,\\
&\Phi\pm\nabla_\beta(u-\widetilde{\varphi})\geq0 &&\text{on }\partial M_\delta,
\end{aligned}
\right.
\end{equation*}
for suitably chosen positive constants $A_{1}$, $A_{2}$, $A_{3}$, where $\mathcal{L}$ and $v$ are defined in Lemma \ref{role-lemma}.
Fix $x_0\in\partial M$ and, for $\delta>0$ sufficiently small, choose a smooth local orthonormal frame $e_1,\ldots,e_n$ on $\overline{M_\delta}$ such that $e_n=\nabla d$ and $e_1,\ldots,e_{n-1}$ are tangent to the level
sets of $d$. In particular, $e_n$ is the inward unit normal on $\partial M\cap B_\delta(x_0)$. By the definition of $\rho$, we have
\begin{equation}\label{5-eq1.5}
\mathcal L(\rho^{2})=T^{ij}\nabla_{ij}(\rho^{2})-f_{p_{i}}\nabla_{i}(\rho^{2})\leq C(1+\sum_{i}T^{ii}).
\end{equation}
By \eqref{3-eq7-1}, \eqref{Eulereq-2}, and Lemma \ref{lemma3}, we obtain
\begin{equation}\label{5-eq2}
\begin{aligned}
\left|\mathcal L\bigl(\nabla_\beta(u-\widetilde{\varphi})\bigr)\right|&\leq~2|T^{ij}\Gamma^l_{i\beta}u_{jl}|+C\big(1+\sum_{i=1}^{n}F^{ii}\big) \\
&\leq C\left(1+\sum_{i=1}^{n}\frac{\partial\widehat T}{\partial\lambda_i}+\sum_{i=1}^{n}\frac{\partial\widehat T}{\partial\lambda_i}\,|\lambda_i|\right),
\end{aligned}
\end{equation}
where $\lambda_1,\ldots,\lambda_n$ are the eigenvalues of $\nabla^2u$, and $\Gamma^l_{ij}=\langle\nabla_{e_i}e_j,e_l\rangle$ are the connection coefficients of the chosen frame.
For each $\gamma<n$, we have
\begin{equation}\label{5-eq3}
\begin{aligned}
\mathcal L\bigl(|\nabla_\gamma(u-\widetilde{\varphi})|^2\bigr)
=&~2T^{ij}\nabla_\gamma(u-\widetilde{\varphi})\nabla_{ij}\nabla_\gamma(u-\widetilde{\varphi})
+2T^{ij}\nabla_i\nabla_\gamma(u-\widetilde{\varphi})\nabla_j\nabla_\gamma(u-\widetilde{\varphi})\\
&-2f_{p_i}\nabla_\gamma(u-\widetilde{\varphi})\nabla_i\nabla_\gamma(u-\widetilde{\varphi})\\
\geq&~2T^{ij}u_{i\gamma}u_{j\gamma}-C\left(1+\sum_{i=1}^{n}\frac{\partial\widehat T}{\partial\lambda_i}+\sum_{i=1}^{n}\frac{\partial\widehat T}{\partial\lambda_i}|\lambda_i|\right).
\end{aligned}
\end{equation}
By \cite[Proposition 2.19]{Guan14}, we know that there exists an index $r$ such that
\begin{equation}\label{5-eq4}
\sum_{\gamma<n}T^{ij}u_{i\gamma}u_{j\gamma}\ge\frac{1}{2}\sum_{i\ne r}\frac{\partial\widehat T}{\partial\lambda_i}\lambda_i^{2}.
\end{equation}
Moreover,
\begin{equation*}
\frac{\partial\widehat T}{\partial\lambda_i}>0,\qquad\sum_i\frac{\partial\widehat T}{\partial\lambda_i}=\sum_iT^{ii}.
\end{equation*}
By \eqref{Eulereq-1} and the estimate for $\tilde F_\alpha$ established in the proof of Lemma~\ref{lemma2}, we have
\begin{equation*}
\sum_i\frac{\partial\widehat T}{\partial\lambda_i}\lambda_i=T^{ij}u_{ij}=f-\alpha\tilde F_\alpha,\qquad\left|\sum_i\frac{\partial\widehat T}{\partial\lambda_i}\lambda_i\right|\leq C_{0}\left(1+\sum_iT^{ii}\right).
\end{equation*}
Consequently, for index $r$ and any $\varepsilon>0$, we have
\begin{equation}\label{5-eq5}
\begin{aligned}
\sum_i\frac{\partial\widehat T}{\partial\lambda_i}|\lambda_i|&\leq 2\sum_{i\neq r}\frac{\partial\widehat T}{\partial\lambda_i}|\lambda_i|+\left|\sum_i\frac{\partial\widehat T}{\partial\lambda_i}\lambda_i\right| \\
&\leq\varepsilon\sum_{i\neq r}\frac{\partial\widehat T}{\partial\lambda_i}\lambda_i^2+\frac{1}{\varepsilon}\sum_{i\neq r}\frac{\partial\widehat T}{\partial\lambda_i}+\left|\sum_i\frac{\partial\widehat T}{\partial\lambda_i}\lambda_i\right|\\
&\leq\varepsilon\sum_{i\neq r}\frac{\partial\widehat T}{\partial\lambda_i}\lambda_i^2+(C_0+\frac{1}{\varepsilon})\left(1+\sum_iT^{ii}\right),
\end{aligned}
\end{equation}
where we used Young's inequality in the second inequality.
Combining \eqref{5-eq1.5}-\eqref{5-eq5},
we obtain
\begin{equation*}
\begin{aligned}
&A_2\mathcal L(\rho^2)-A_3\mathcal L\big(\sum_{\gamma<n}|\nabla_\gamma(u-\widetilde{\varphi})|^2\big)\pm\mathcal L\bigl(\nabla_\beta(u-\widetilde{\varphi})\bigr)\\
\leq&~C_{1}A_2\left(1+\sum_iT^{ii}\right)+C_{1}\left(1+\sum_iT^{ii}+\sum_i\frac{\partial\widehat T}{\partial\lambda_i}|\lambda_i|\right)\\
&-2A_{3}\sum_{\gamma<n}T^{ij}u_{i\gamma}u_{j\gamma}+C_{1}A_{3}(n-1)\left(1+\sum_{i=1}^{n}\frac{\partial\widehat T}{\partial\lambda_i}+\sum_{i=1}^{n}\frac{\partial\widehat T}{\partial\lambda_i}|\lambda_i|\right)\\
\leq&~ C_{1}A_2\left(1+\sum_iT^{ii}\right)-A_3\sum_{i\neq r}\frac{\partial\widehat T}{\partial\lambda_i}\lambda_i^2
+C_{1}\big((n-1)A_{3}+1\big)\left(1+\sum_iT^{ii}+\sum_i\frac{\partial\widehat T}{\partial\lambda_i}|\lambda_i|\right)\\
\leq&~\bigl(C_{1}(n-1)\varepsilon A_{3}+C_{1}\varepsilon-A_3\bigr)\sum_{i\neq r}\frac{\partial\widehat T}{\partial\lambda_i}\lambda_i^2+\bigl(C_{1}A_2+(C_{2}+\frac{C_{2}}{\varepsilon})\big(A_{3}+1\big)\bigr)\left(1+\sum_iT^{ii}\right),
\end{aligned}
\end{equation*}
where $C_{2}$ is a positive constant depending only on $n$, $C_{0}$ and $C_{1}$. By choosing $0<\varepsilon<\min\{\frac{1}{2C_{1}(n-1)}, 1\}$ and $A_{3}>2C_{1}\varepsilon$, we have
\begin{equation*}
A_2\mathcal L(\rho^2)-A_3\mathcal L\big(\sum_{\gamma<n}|\nabla_\gamma(u-\widetilde{\varphi})|^2\big)\pm\mathcal L\bigl(\nabla_\beta(u-\widetilde{\varphi})\bigr)
\leq\bigl(C_{1}A_2+C_{3}A_{3}\bigr)\left(1+\sum_iT^{ii}\right),
\end{equation*}
where $C_{3}$ is a positive constant depending only on $n$, $C_{1}$ and $C_{2}$. Combining the above inequality with Lemma~\ref{role-lemma} and choosing $A_1\gg A_2\gg A_3\gg1$, we complete the proof
of the claim.

Then the maximum principle gives
\begin{equation*}
\Phi\pm \nabla_{\beta}(u-\widetilde{\varphi})\ge0\qquad\text{in } M_{\delta}.
\end{equation*}
Taking their inward normal derivatives at $x_{0}$, for $\beta=1,\dots, n-1$, we obtain
\begin{equation*}
|\nabla_{n\beta}u(x_{0})|\leq \nabla_{n}\Phi(x_{0})+|\nabla_{n\beta}\widetilde{\varphi}(x_{0})|\leq C.
\end{equation*}

\textbf{Case 3}: Estimate of $\nabla_{nn}u$ on $\partial M$.

Since $\eta[U]\in\Gamma_{k-1}$ and $k\ge2$, we have $(n\tau-1)\Delta u=\sigma_{1}(\eta[U])>0$ in $M$. Choose the boundary frame $\{e_{1},\dots,e_{n}\}$ as in Case 2. Then by
continuity and the tangential estimate from Case 1, we have
\begin{equation*}
u_{nn}\ge-\sum_{\gamma<n}u_{\gamma\gamma}\ge-C\qquad\text{on } \partial M.
\end{equation*}
It remains to prove an upper bound for $u_{nn}$.

On $\partial M$, the matrix of $U$ has the block form
\begin{equation*}
U=
\begin{pmatrix}
\bigl[\tau(u_{nn}+\sum_{\iota<n}u_{\iota\iota})\delta_{\beta\gamma}
 -u_{\beta\gamma}\bigr]_{\beta,\gamma<n}
 & [-u_{\beta n}]_{\beta<n}\\[2pt]
 [-u_{n\gamma}]_{\gamma<n}
 & (\tau-1)u_{nn}+\tau\sum_{\iota<n}u_{\iota\iota}
\end{pmatrix}.
\end{equation*}
By Cases 1 and 2, $u_{\beta\gamma}$ and $u_{\beta n}$, $\beta,\gamma<n$, are uniformly bounded on $\partial M$. For sufficiently large $u_{nn}$ (otherwise, the desired upper
bound already holds), expanding the principal minors, we obtain
\begin{equation}\label{C3-expansion}
S_{j}(U)=\left[C_{n-1}^{j}\tau^{j}+C_{n-1}^{j-1}\tau^{j-1}(\tau-1)\right]u_{nn}^{j}+O(u_{nn}^{j-1}),\qquad 1\le j\le n.
\end{equation}
where the constant implicit in the error term is uniform on $\partial M$, $C_{n-1}^{n}=0$ for $j=n$, and $S_{0}=1$ for $j=0$.

If $k<n$, the coefficient of $u_{nn}^{k}$ in \eqref{C3-expansion} is positive for every $\tau\ge1$. The same holds if $k=n$ and $\tau>1$. In either case, for $u_{nn}$ sufficiently large, we have
\begin{equation*}
f^{k-l}=\frac{S_{k}(U)}{S_{l}(U)}\ge c\,u_{nn}^{k-l}.
\end{equation*}
Since $f(x,u,\nabla u)$ is bounded by Lemmas~\ref{lemma1} and~\ref{lemma3}, the desired upper bound follows.

It remains to consider the case $k=n$ and $\tau=1$, in which the coefficient of $u_{nn}^{n}$ in the expansion of $S_n(U)$ vanishes.
Expanding $S_n(U)=\sigma_{n}(U)+\alpha\sigma_{n-1}(U)$, we obtain
\begin{equation}\label{C3-critical}
S_n(U)=\left(\sum_{\gamma<n}u_{\gamma\gamma}+\alpha\right)u_{nn}^{n-1}+O(u_{nn}^{n-2}).
\end{equation}
We claim that
\begin{equation}\label{C3-trace-lower}
m_{0}:=\min_{\partial M}\Big(\sum_{\gamma<n}u_{\gamma\gamma}+\alpha\Big)\ge c_{0}>0
\end{equation}
with $c_{0}$ independent of the second derivatives of $u$.

Since $\eta[U]\in\widetilde{\Gamma}_n$ is equivalent to $(\eta[U],\alpha)\in\Gamma_n^{n+1}$, the properties of Garding cones imply that deleting $n-1$ coordinates gives $\eta_i+\alpha>0$ for every $i$, so $[U+\alpha g]$ is positive definite. The Newton--Maclaurin inequalities \eqref{prop-5-ieq2} give
\begin{equation}\label{6-eq4}
S_j(U)\ge c(n,l,j)f(x,u,\nabla u)^j,\qquad 1\le j\le n,
\end{equation}
where $c(n,l,j)$ is a positive constant depending only on $n$, $l$, $j$, while the $C^0$ and $C^1$ estimates and the positivity of $f$ give a uniform positive lower bound for $f(x,u,\nabla u)$.
For each fixed solution, \eqref{6-eq4} extends to $\partial M$ by continuity, so $(\eta[U],\alpha)\in\Gamma_n^{n+1}$ and $[U+\alpha g]$ is positive definite on $\partial M$. Thus we obtain a positive lower bound $m_0$ for each fixed solution.

We next prove that this positive lower bound can be chosen independently of the solution.
For any solution $u$, let $x_1\in\partial M$ be a point where $m_{0}:=\min_{\partial M}\bigl(\sum_{\gamma<n}u_{\gamma\gamma}+\alpha\bigr)$ is attained. Since the subsolution $\underline{u}$ is $n$-admissible, there exists a constant $\kappa_0>0$ such that
\begin{equation*}
\underline{U}_{nn}+\alpha\ge\kappa_{0}>0 \qquad\text{on}~\partial M.
\end{equation*}
If $m_{0}\ge\kappa_{0}/2$, the claim follows. Suppose $m_{0}<\kappa_{0}/2$, let
\begin{equation*}
B_{\beta\gamma}=\langle\nabla_{e_{\beta}}e_{\gamma},e_{n}\rangle,\qquad H=\sum_{\gamma<n}B_{\gamma\gamma}.
\end{equation*}
Since $u=\underline{u}$ on $\partial M$, we have $(u-\underline u)_{\gamma\gamma}=-(\nabla_{e_\gamma}e_\gamma)(u-\underline u)$, and hence
\begin{equation*}
\nabla_{n}(u-\underline{u})(x_{1}) H(x_{1})=\underline{U}_{nn}(x_{1})+\alpha-m_{0}\ge\frac{\kappa_{0}}{2}.
\end{equation*}
Let $\Omega_{\delta_0}$ and $M_\delta$ be as in the proof of Lemma~\ref{role-lemma}, with $M_\delta$ centered at $x_1$. By Lemma \ref{lemma1} and Lemma \ref{lemma3} we have $0\le\nabla_{n}(u-\underline{u})\le C$ on $\partial M$, hence $H(x_{1})\ge c_{0}>0$. By the smoothness of $\partial M$, $H\ge c_0/2$
on $\partial M\cap B_\delta(x_1)$ for sufficiently small $\delta>0$. On $\partial M\cap B_\delta(x_1)$, let
\begin{equation*}
\Theta=\frac{\Delta_{\partial M}\varphi+\alpha-m_0}{H}.
\end{equation*}
For sufficiently small $\delta_0>0$, extend $\Theta$ to the corresponding part of $\Omega_{\delta_0}$ by keeping it constant along the normal geodesics, so that $\nabla_n\Theta=0$.
Since $\varphi\in C^4(\partial M)$,
$0<m_0<\kappa_0/2$, and $H\ge c_0/2$, the $C^2$ norm of $\Theta$ is uniformly bounded on $\overline{M_\delta}$ for sufficiently small $\delta>0$.
By $\sum_{\gamma<n}u_{\gamma\gamma}=\Delta_{\partial M}\varphi-Hu_n$ on $\partial M$, we obtain
\begin{equation*}
\Theta-u_n=\frac{\sum_{\gamma<n}u_{\gamma\gamma}+\alpha-m_0}{H}\ge0\qquad\text{on}~\partial M\cap B_\delta(x_1),
\end{equation*}
with equality at $x_1$. Then we use the auxiliary function $\Phi$ defined in \eqref{5-eq1}, now centered at $x_1$, with $\rho(x)=\operatorname{dist}(x,x_1)$ and $M_\delta=M\cap B_\delta(x_1)$.
As in \eqref{5-eq2}, we obtain
\begin{equation*}
\begin{aligned}
|\mathcal L(\Theta-u_n)|&\leq 2|T^{ij}\Gamma^l_{in}u_{jl}|+C\left(1+\sum_iT^{ii}\right)\\
&\leq C\left(1+\sum_iT^{ii}+\sum_i\frac{\partial\widehat T}{\partial\lambda_i}|\lambda_i|\right).
\end{aligned}
\end{equation*}
Combining this with \eqref{5-eq1.5} and \eqref{5-eq3}--\eqref{5-eq5}, we obtain
\begin{equation*}
\begin{aligned}
\mathcal L(\Phi+\Theta-u_n)&\leq A_1\mathcal Lv-\bigl(A_3-C_1((n-1)A_3+1)\varepsilon\bigr)\sum_{i\ne r}\frac{\partial\widehat T}{\partial\lambda_i}\lambda_i^2\\
&\quad+\left(C_1A_2+\left(C_2+\frac{C_2}{\varepsilon}\right)(A_3+1)\right)\left(1+\sum_iT^{ii}\right).
\end{aligned}
\end{equation*}
As in Case 2, choosing $\varepsilon>0$ sufficiently small and $A_1\gg A_2\gg A_3\gg1$, and applying Lemma~\ref{role-lemma}, we obtain
\begin{equation*}
\mathcal L(\Phi+\Theta-u_n)\le0\qquad\text{in }M_\delta.
\end{equation*}
On $\partial M\cap B_\delta(x_1)$, we have $v=0$, $\nabla_\gamma(u-\widetilde{\varphi})=0$ for $\gamma<n$,
and $\Theta-u_n\ge0$, so $\Phi+\Theta-u_n\ge0$.
On $M\cap\partial B_\delta(x_1)$, for sufficiently large $A_2$, we have $\Phi+\Theta-u_n\ge A_2\delta^2-C(A_3+1)\ge0$. Thus,
\begin{equation*}
\left\{
\begin{aligned}
&\mathcal L(\Phi+\Theta-u_n)\le0
&&\text{in }M_\delta,\\
&\Phi+\Theta-u_n\ge0
&&\text{on }\partial M_\delta.
\end{aligned}
\right.
\end{equation*}
The maximum principle yields $\Phi+\Theta-u_n\ge0$ in $M_\delta$. Since this function vanishes at $x_1$, $\nabla_n\Theta=0$, and $\nabla_{e_n}e_n=0$, taking its inward normal derivative at $x_1$ gives
\begin{equation*}
u_{nn}(x_1)\le\nabla_n\Phi(x_1)=A_1\bigl(\nabla_n(u-\underline u)(x_1)+t\bigr)\le C.
\end{equation*}
By Cases 1 and 2 and the bounds for $u_{nn}(x_1)$, all components of $U(x_1)$ are uniformly bounded. Together with \eqref{6-eq4} and the uniform positive lower bound for $f$, this implies that $(\eta[U(x_1)],\alpha)$ lies in a
fixed compact subset of $\Gamma_n^{n+1}$. Hence
\begin{equation*}
m_0=U_{nn}(x_1)+\alpha\ge\min_i\bigl(\eta_i[U(x_1)]+\alpha\bigr)\ge c_0>0,
\end{equation*}
where $c_0$ is independent of the solution. This proves \eqref{C3-trace-lower}.

Finally, by \eqref{C3-critical} and \eqref{C3-trace-lower}, for sufficiently large $u_{nn}$ we have
\begin{equation*}
S_n(U)\ge \frac{c_{0}}{2}\,u_{nn}^{n-1},\qquad0<S_l(U)\le C_{0} u_{nn}^l,
\end{equation*}
where $S_0=1$ when $l=0$. Hence
\begin{equation*}
f^{n-l}=\frac{S_n(U)}{S_l(U)}\ge c\,u_{nn}^{n-1-l}.
\end{equation*}
Since $l\le n-2$ and $f(x,u,\nabla u)$ is bounded by Lemmas~\ref{lemma1} and~\ref{lemma3}, this gives the uniform upper bound for $u_{nn}$ on $\partial M$. Together with Cases 1 and 2, we conclude that
\begin{equation*}
\sup_{\partial M}|\nabla^2u|\le C.
\end{equation*}
The proof is complete.
\end{proof}

\begin{proof}[Proof of Theorem \ref{maintheorem1} and Theorem \ref{maintheorem2}]
The conclusions follow from Lemma \ref{lemma1}, Lemma \ref{lemma3}  and Lemma \ref{lemma5}.
\end{proof}

\begin{proof}[Proof of Theorems~\ref{maintheorem3} and~\ref{maintheorem4}]
Theorems~\ref{maintheorem1} and~\ref{maintheorem2} give uniform $C^2(\overline M)$ estimates for classical elliptic solutions of the Dirichlet problems:
\begin{equation*}
\left\{
\begin{aligned}
&F(U[u])=t f(x,u,\nabla u)+(1-t)F(U[\underline u])
&&\text{in }M,\\
&u=\varphi &&\text{on }\partial M,
\end{aligned}
\right.
\end{equation*}
for $0\le t\le1$. By the Newton--Maclaurin inequalities, the uniform $C^2$ bounds and the positive lower bound for the right-hand side imply that $(\eta[U],\alpha)$ remains in a fixed compact subset of $\Gamma_k^{n+1}$, so the equation is uniformly elliptic. Since $u=\underline u$ solves the problem at $t=0$, existence follows by the method of continuity, the Evans--Krylov \cite{EK82, Kry82} second derivative H\"older estimates, and the boundary regularity
estimates, with a standard approximation argument for a $C^2$ subsolution; see \cite{CNS85,CTX23,GT01}. Uniqueness follows from the maximum principle.
\end{proof}

\begin{remark}
Recently, some researchers have begun to study geometric and PDE
problems for equations of the following general sum form:
\begin{equation*}
F(\lambda)=\sigma_{k}^{(n)}(\lambda)+\sum_{r=1}^{m}a_{r}\sigma_{k-r}^{(n)}(\lambda),
\end{equation*}
where $\sigma_{k}^{(n)}$ denotes the $k$-th elementary symmetric
function on $\mathbb{R}^{n}$. A natural question is whether
problems of the following form can be considered on Riemannian
manifolds:
\begin{equation*}
\left\{
\begin{aligned}
&\sigma_{k}^{(n)}(\eta[U])+\sum_{r=1}^{m}a_{r}\sigma_{k-r}^{(n)}(\eta[U])=f(x,u,\nabla u),\qquad&&~\mathrm{in}~M,\\
&u=\varphi, \qquad&&~\mathrm{on}~ \partial M.
\end{aligned}
\right.
\end{equation*}
This is the Dirichlet boundary value problem that we will study
in an upcoming paper.
\end{remark}

\end{document}